\documentclass[11pt,twoside]{amsart}
\usepackage{amsmath,amssymb,latexsym}

\newtheorem{theorem}{Theorem}[section]
\newtheorem{lemma}[theorem]{Lemma}
\newtheorem{definition}[theorem]{Definition}
\newtheorem{corollary}[theorem]{Corollary}

\newtheorem{construction}[theorem]{Construction}
\renewcommand{\epsilon}{\varepsilon}

\newcommand{\N}{\mathbb{N}}

\begin{document}

\thispagestyle{empty}

\title{Extended visual cryptography schemes}
\author{Andreas Klein and Markus Wessler}

\begin{abstract}
Visual cryptography schemes have been introduced in 1994 by Naor and
Shamir. Their idea was to encode a secret image into $n$ shadow
images and to give exactly one such shadow image to each member of a
group $P$ of $n$ persons. Whereas most work in recent years has been done
concerning the problem of qualified and forbidden subsets of $P$ or
the question of contrast optimizing, in this paper we study
extended visual cryptography schemes, i.e. shared secret
systems where any subset of $P$ shares its own secret.
\end{abstract}

\maketitle


\section{\sc Introduction}
A visual cryptography scheme is given by the following set up. Let $P$
be a group of $n$ persons where each participant is given exactly
one image (in fact it does not have to be a real image) xeroxed onto a
transparency. Stacking all the transparencies
together, a secret image is recovered. So in this sense the
participants share a secret. This set up can be generalized to the
case where some subsets $X \subseteq P$ (which are usually called
{\em qualified subsets of $P$}) can recover the secret by stacking their
transparencies together, whereas other, {\em forbidden} subsets
cannot. Such structures, called {\em access structures}, have been
examined very well.
In \cite{Naor:1994} Naor and Shamir analysed so-called $(k,n)$-threshold
visual cryptography schemes, i.e. schemes where a subset is qualified
if and only if it consists of at least $k$ participants. In
\cite{Ateniese:1996II} and \cite{Ateniese:1996I}
their idea was extended to general access structures.

Most work concerning this subject focuses on two aspects, either the
pixel expansion, i.e. the number of subpixels which is needed on the
different levels to represent a white or a black pixel,
or the contrast, i.e. the difference
of subpixels representing a white or a black pixel.

As a further generalization, the existence of a secret
image can be concealed by displaying a different image on each
slide. Naor and Shamir \cite{Naor:1994} solved this problem for the
$(2,2)$-threshold scheme. In \cite{Ateniese:2001} this problem
was considered for a general access structure.
In \cite{Droste:1996} Droste made
a further generalization: Stacking
the transparencies of each participant together, a secret image is
recovered, and there is in fact only this single way to recover it. But
moreover, the participants of any arbitrary subset $X$ of $P$ share a
secret, too. Hence we have $2^n - 1$ more or less secret images.

We start by briefly recalling the work done by
Droste and prove that the scheme proposed in \cite{Droste:1996} has
minimal pixel expansion. Then we prove a trade-off theorem between the
contrast of the different images.

Finally we give new constructions for generalized visual cryptography
schemes with less then $2^n - 1$ subsets in order to achieve a smaller
pixel expansion and a better contrast.


\section{Preliminaries}
A visual cryptography scheme is based on the fact that each pixel of
an image is divided into a certain number $m$ of subpixels. This
number $m$ is called the {\em pixel expansion} of the image. If the
number of black subpixels needed to represent a white pixel in an
image is $l$, and the number of black subpixels needed to represent a
black pixel is $h$, then we call the number $\alpha = \frac{h-l}{m}$
the {\em contrast} of the image.

An extended visual cryptography scheme consists of $n$
transparencies $\tau_1 , \dots , \tau_n$ and $2^{n}-1$ different
images (one for each non-empty subset
$T \subseteq \{1 , \cdots , n\}$). We denote by $I_T$ the image
which is recovered by stacking together exactly the transparencies
$\tau_i$ for $i \in T$. We generalize this as follows.
For any non-empty subset $\mathfrak{S} \subseteq \mathcal{P}(\{1 ,
\cdots , n\})
\backslash \{\emptyset\}$
an $\mathfrak{S}$-extended visual cryptography scheme consists of $n$
transparencies $\tau_1 , \dots , \tau_n$ with the following property:
Let $T \in \mathfrak{S}$. If we stack together the slides $\tau_i$ for
$i \in T$,
then we recover the image $I_T$ for which each white pixel is
represented by $l_T$ black subpixels and each black pixel is
represented by $h_T$ black subpixels.
Furthermore, for $T'$ not contained in $T$, the distribution of
subpixels on the transparencies $\tau_i$ with $i \in T$ is independent
of the image $I_{T'}$, i.e. the information of the transparencies
$\tau_i$ with $i \in T$ does not suffice to recover the image
$I_{T'}$.

More formally, we define an $\mathfrak{S}$-extended visual cryptography
scheme
as follows. (See also \cite{Naor:1994} for ``usual'' visual cryptography
schemes and \cite{Droste:1996} for $\mathfrak{S}$-extended visual
cryptography schemes.)

\begin{definition}
  Let $\mathfrak{S} \subseteq \mathcal{P}(\{1,\dots,n\}) \backslash
  \{\emptyset\}$.

  An $\mathfrak{S}$-extended visual cryptography scheme is described by
  multi-sets
  $C^{\mathfrak{T}}$ of $n \times m$ Boolean matrices for
  $\mathfrak{T} \subseteq \mathfrak{S}$. (For given $\mathfrak{T}$
  each Boolean matrix in $C^{\mathfrak{T}}$
  describes the colors of the subpixels on each transparency, where the
  corresponding pixel in image $I_T$ is black if and only if $T \in
  \mathfrak{T}$. For encoding, each matrix in $C^{\mathfrak{T}}$ is
  chosen with the same probability.)

  The
  multi-sets $C^{\mathfrak{T}}$ must satisfy the following conditions:
  \begin{enumerate}
  \item Let $B \in C^{\mathfrak{T}}$. For $\{i_{1},\dots,i_{q}\} \in
    \mathfrak{S}$ the
    Hamming weight of the OR of the rows $i_{1},\dots,i_{q}$ of $B$ is
    $h_{\{i_{1},\dots,i_{q}\}}$ if $\{i_{1},\dots,i_{q}\} \in
\mathfrak{T}$ and
    $l_{\{i_{1},\dots,i_{q}\}}$ otherwise, i.e.
    $$ w_{Ham}((b_{i_{1},1},\dots,b_{i_{1},m})\text{ OR }\dots
    \text{ OR }(b_{i_{q},1},\dots,b_{i_{q},m})) =
    \begin{cases}
      h_{\{i_{1},\dots,i_{q}\}} & \text{if } \{i_{1},\dots,i_{q}\} \in
      \mathfrak{T} \\
      l_{\{i_{1},\dots,i_{q}\}} & \text{if } \{i_{1},\dots,i_{q}\}
      \notin \mathfrak{T}
    \end{cases} \quad.
    $$
    (This means stacking the transparencies
    $\tau_{i_{1}},\dots,\tau_{i_{q}}$ together we recover the image
    $I_{\{i_{1},\dots,i_{q}\}}$.) \smallskip
  \item For $\{i_{1},\dots,i_{q}\} \subseteq \{1,\dots,n\}$ and
    $\mathfrak{T},\mathfrak{T}'
    \subseteq \mathfrak{S}$ with $\mathfrak{T} \cap
    \mathcal{P}(\{i_{1},\dots,i_{q}\}) = \mathfrak{T}' \cap
    \mathcal{P}(\{i_{1},\dots,i_{q}\})$ we obtain the same multi-sets if
we
    restrict the matrices in $C^{\mathfrak{T}}$ and
    $C^{\mathfrak{T}'}$ respectively to the rows
    $i_{1},\dots,i_{q}$.

    (This condition guarantees the security of the different images.)
  \end{enumerate}
If $\mathfrak{S} = \mathcal{P}(\{1, \dots, n\}) \backslash
\{\emptyset\}$ we simply call this an extended visual cryptography
scheme.
\end{definition}


In \cite{Droste:1996} Droste gives the following construction for
$\mathfrak{S}$-extended visual cryptography schemes using
$(k,k)$-threshold schemes.

\begin{construction}\label{droste}
For each $T \in \mathfrak{S}$ we take $2^{|T|-1}$ subpixels and use
them to construct a $(|T|, |T|)$-threshold visual cryptography
scheme. If $i \notin T$ the corresponding subpixels on $\tau_{i}$ will
be black. The $\mathfrak{S}$-extended visual cryptography scheme is
achieved by putting all these schemes together. Since we shall not need
the details of this construction in the sequel, we omit a formal
definition and refer to \cite{Droste:1996}.
\end{construction}

The scheme obtained by this construction has pixel expansion
$$ m = \sum_{T \in \mathfrak{S}} 2^{|T|-1}$$
and the contrast of all encoded images is $\frac{1}{m}$. We shall prove
in the following sections that this construction is optimal if
$\mathfrak{S} = \mathcal{P}(\{1,\dots,n\}) \backslash \{\emptyset\}$
but it is not optimal for general $\mathfrak{S}$.


\section{Pixel Expansion and Contrast for the Extended Scheme}
It is sufficient to consider the case of one single pixel.
For a given non-empty subset $T \subseteq \{1 , \ldots ,
n\}$ let $x_T$ be the number of
subpixels which are black
exactly on the transparencies $i$ for $i \in T$ and let us denote by
$x$ the vector of all the $x_T$.
For ${\emptyset \neq S \subseteq \{1 , \ldots ,
n\}}$ let $r_S$ be the number of black subpixels needed for the
image $I_S$. Formally we set $r_{\emptyset}=0$.
We write $r$ for the vector of all the $r_S$ with $\emptyset \neq S
\subseteq \{1 , \ldots , n\}$.

This leads to the linear equation system given by

\begin{equation}\label{matrix}
   M x = r
\end{equation}

where $M=(m_{S,T})_{\emptyset \neq S,T \subseteq \{1,\dots,n\}}$ is
defined by $m_{S,T}=1$ if $S \cap T \neq \emptyset$ and $m_{S,T}=0$
otherwise.

\begin{lemma} \label{inv-m}
  Equation (\ref{matrix}) has a unique integral solution.
\end{lemma}

\begin{proof}
We prove this by induction on the number $n$ of transparencies.
Writing $M_1 = (1)$, we obtain the following recursion formula which
follows directly from the definition:
$$ M_{n+1} =
\begin{pmatrix}
  M_{n} & {\bf 0} & M_{n} \\
  {\bf 0} & 1 & {\bf 1} \\
  M_{n} & {\bf 1} & {\bf 1}
\end{pmatrix}.
$$
Here the index denotes the number of transparencies.

Let $e_{n}$ denote the $(2^{n}-1)$ dimensional vector
$(1,\dots,1)^{t}$. With $M_{1}^{-1} = (1)$
we obtain the following recursion formula for $M_{n}^{-1}$:
$$ M_{n+1}^{-1} =
\begin{pmatrix}
  0 & -M_{n}^{-1}e_n & M_{n}^{-1} \\
  -e_n^{t}M_{n}^{-1} & 0 & e_n^{t}M_{n}^{-1} \\
  M_{n}^{-1} & M_{n}^{-1}e_n & -M_{n}^{-1}
\end{pmatrix}
$$
We notice that the components of
$M_n^{-1}e_n$ are only $-1, 0$ and $1$ and that
$e_{n}^{t}M_{n}^{-1}e_{n}=1$. Then the formula can be proved by induction.

Thus $M_{n}^{-1}$ contains only the entries $-1,0$ and $1$ and
therefore equation (\ref{matrix}) has an integral solution.
\end{proof}

\begin{lemma}
  The solution of (\ref{matrix}) is non-negative if and only if for
  each $S \subsetneq \{1,\dots,n\}$ the condition
  \begin{equation}
    \label{eq:cond}
    \sum_{S \subseteq T \subseteq \{1,\dots,n\}} (-1)^{|S|+|T|} r_{T}
    \le 0
  \end{equation}
  is satisfied.
\end{lemma}
\begin{proof}
  We claim that $x = (x_{S})$ with
  \begin{equation}
    \label{eq:X-sol}
    x_{S} = \sum_{\{1,\dots,n\}\backslash S \subseteq T \subseteq
      \{1,\dots,n\}} (-1)^{|T|+|S|+n+1} r_{T}
  \end{equation}
  solves equation (\ref{matrix}) and due to Lemma
  \ref{inv-m} this solution is unique.

  To prove this we substitute $x$ in equation (\ref{matrix}). For
  $\emptyset \neq U \subseteq \{1,\dots,n\}$ the line of the system
  of linear equations corresponding to $U$ yields
  \begin{equation}
    \label{matrix2}
    \begin{split}
      \sum_{\emptyset \neq S \subseteq \{1,\dots,n\}} m_{U,S}x_{S} &=
      \sum_{S \subseteq \{1,\dots,n\}}
      m_{U,S}\sum_{\{1,\dots,n\}\backslash S \subseteq T \subseteq
        \{1,\dots,n\}} (-1)^{|T|+|S|+n+1} r_{T} \\ \\
      &=\sum_{T \subseteq \{1,\dots,n\}} \quad\sum_{\{1,\dots,n\}
\backslash T
      \subseteq S \subseteq \{1,\dots,n\}} (-1)^{|T|+|S|+n+1} r_{T}
      m_{U,S} \\ \\
      &=\sum_{\emptyset \neq T \subseteq \{1,\dots,n\}} (-1)^{|T|}r_{T}
      \sum_{\{1,\dots,n\} \backslash T
      \subseteq S \subseteq \{1,\dots,n\}} (-1)^{|S|+n+1} m_{U,S} \quad.
    \end{split}
  \end{equation}

  If $T \not\subseteq U$ we choose $t \in T \backslash U$ and obtain
  \begin{multline*}
      \sum_{\{1,\dots,n\} \backslash T \subseteq S \subseteq
        \{1,\dots,n\}}
      (-1)^{|S|+n+1}m_{U,S}  \\
      = \sum_{\{1,\dots,n\} \backslash \{T \cup \{t\}\} \subseteq S
        \subseteq \{1,\dots,n\} \backslash \{t\}}
      (-1)^{|S|+n+1}  m_{U,S} +  (-1)^{|S|+1+n+1}m_{U,S \cup \{t\}}
      = 0
  \end{multline*}
  since $m_{U,S} = m_{U,S \cup \{t\}}$.

  If $T \subsetneq U$ and $T \neq \emptyset$ we find
  \begin{equation*}
    \begin{split}
      \sum_{\{1,\dots,n\} \backslash T \subseteq S \subseteq
\{1,\dots,n\}}  (-1)^{|S|+n+1}m_{U,S}
      &= \sum_{\{1,\dots,n\} \backslash T \subseteq S \subseteq
\{1,\dots,n\}} (-1)^{|S|+n+1}
      = 0
    \end{split} \quad
  \end{equation*}
  since $m_{U,S} = 1$.

  But for $\emptyset \neq T=U$ we find
  \begin{equation*}
    \begin{split}
      \sum_{\{1,\dots,n\} \backslash T \subseteq S \subseteq
\{1,\dots,n\}}  (-1)^{|S|+n+1}m_{U,S}
      &= \sum_{\{1,\dots,n\} \backslash T \subseteq S \subseteq
\{1,\dots,n\}} (-1)^{|S|+n+1}
      - (-1)^{|\{1,\dots,n\}\backslash U|+n+1} \\
      &= (-1)^{|\{1,\dots,n\}\backslash U|+n}
    \end{split}
  \end{equation*}
  since $m_{U,S} = 1$ for $S \neq \{1,\dots,n\} \backslash U$.

  Thus equation (\ref{matrix2}) yields
  \begin{equation*}
    \begin{split}
      \sum_{\emptyset \neq S \subseteq \{1,\dots,n\}} m_{U,S}x_{S} &=
      \sum_{\emptyset \neq T \subseteq \{1,\dots,n\}} (-1)^{|T|}r_{T}
      \sum_{\{1,\dots,n\} \backslash T \subseteq S \subseteq
      \{1,\dots,n\}}  (-1)^{|S|+n+1} m_{U,S}\\ \\
      &= (-1)^{|U|}r_{U} (-1)^{|\{1,\dots,n\}\backslash U|+n}\\
      &= (-1)^{2n}r_{U} = r_{U} \quad.
    \end{split}
  \end{equation*}

  This proves that $x$ is a solution of equation (\ref{matrix}) and
  the lemma follows.
\end{proof}

Now we can solve~\eqref{eq:cond} to derive bounds for the pixel
expansion and the contrast.

\begin{theorem}\label{drostegut}
  An extended visual cryptography scheme with $n$ transparencies
  needs at least
  $\frac{1}{2}(3^n-1)$ subpixels. Hence Construction
  \ref{droste} is optimal with respect to the pixel expansion.
\end{theorem}
\begin{proof}
  First we note that the inequality~\eqref{eq:cond} is as strong as
  possible if $r_{T} = h_{T}$ for $|S|+|T|$ even and $r_{T}=l_{T}$ for
  $|S|+|T|$ odd.

  Thus an extended visual cryptography scheme exists if
  and only if
  \begin{equation}
    \label{eq:existence}
    \sum_{\substack{S \subseteq T \subseteq \{1,\dots,n\} \\
        |S| \equiv |T| \mod 2}} h_{T} \le
    \sum_{\substack{S \subseteq T \subseteq \{1,\dots,n\} \\
        |S| \not\equiv |T| \mod 2}} l_{T}
  \end{equation}
  holds for each $S \subsetneq \{1,\dots,n\}$.

  For each $\emptyset \neq T \subseteq \{1,\dots,n\}$ let
  $\delta_{T}=h_{T}-l_{T}$. Our goal is to prove that
  \begin{equation}
    \label{eq:m-bound}
    m \ge h_{\{1,\dots,n\}} \ge \sum_{\emptyset \neq T \subseteq
    \{1,\dots,n\}} \delta_{T}2^{|T|-1} \quad.
  \end{equation}

  As the first step in this proof we show that, for given values
  $\delta_T$ (for $\emptyset \neq T \subseteq \{1, \dots , n\})$, the
number
  $h_{\{1,\dots,n\}}$ is
  minimal if for all $S \subsetneq \{1,\dots,n\}$
  inequality~\eqref{eq:existence} is satisfied with equality.

  To this end, suppose
  $$
  \sum_{\substack{S \subseteq T \subseteq \{1,\dots,n\} \\
      |S| \equiv |T| \mod 2}} h_{T} <
  \sum_{\substack{S \subseteq T \subseteq \{1,\dots,n\} \\
      |S| \not\equiv |T| \mod 2}} l_{T}
  $$
  for some $S \subsetneq \{1,\dots,n\}$. But the contrast levels
  \begin{align*}
    \bar{h}_{T} &=
    \begin{cases}
      h_{T} & \text{for } T \subseteq S \\
      h_{T}-1 & \text{otherwise}
    \end{cases}
    \intertext{and}
    \bar{l}_{T} &=
    \begin{cases}
      l_{T} & \text{for } T \subseteq S \\
      l_{T}-1 & \text{otherwise}
    \end{cases}
  \end{align*}
  satisfy~\eqref{eq:existence}, since
  \begin{multline*}
    |\{T \mid S \subseteq T \subseteq \{1,\dots,n\}; |T|\equiv|S|\mod 2;
      T \not\subseteq S\}| =\\
    |\{T \mid S \subseteq T \subseteq \{1,\dots,n\}; |T|\not\equiv|S|\mod
2;
      T \not\subseteq S\}| \quad.
  \end{multline*}

  Thus we may assume that inequality~\eqref{eq:existence} is satisfied
with
  equality for each  $S \subsetneq \{1,\dots,n\}$.

  Next we claim that
  \begin{equation}
    \label{eq:hT}
      h_{T}
      = \sum_{\emptyset \neq T' \subseteq
        \{1,\dots,n\}}\delta_{T'}2^{|T'|-1}
      - \sum_{T \subsetneq T' \subseteq \{1,\dots,n\}}
      \delta_{T'}2^{|T'|-1-|T|}
  \end{equation}
  for $\emptyset \neq S \subseteq \{1,\dots,n\}$
  satisfy~\eqref{eq:existence} with equality.

  To prove this we have to show that
  \begin{multline*}
    \sum_{\substack{S \subseteq T \subseteq \{1,\dots,n\} \\
        |S| \equiv |T| \mod 2}}
    \left[
      \sum_{\emptyset \neq T' \subseteq
        \{1,\dots,n\}}\delta_{T'}2^{|T'|-1}
      - \sum_{T \subsetneq T' \subseteq \{1,\dots,n\}}
      \delta_{T'}2^{|T'|-1-|T|}
    \right] = \\
    \sum_{\substack{S \subseteq T \subseteq \{1,\dots,n\} \\
        |S| \not\equiv |T| \mod 2}}
    \left(\left[
      \sum_{\emptyset \neq T' \subseteq
        \{1,\dots,n\}}\delta_{T'}2^{|T'|-1}
      - \sum_{T \subsetneq T' \subseteq \{1,\dots,n\}}
      \delta_{T'}2^{|T'|-1-|T|}
    \right] - \delta_{T} \right)
  \end{multline*}
  or equivalently
  \begin{multline*}
    \sum_{\emptyset \neq T' \subseteq \{1,\dots,n\}}\delta_{T'}
    \left[
     \sum_{\substack{S \subseteq T \subseteq \{1,\dots,n\} \\
         |S| \equiv |T| \mod 2}}  2^{|T'|-1}
     - \sum_{\substack{S \subseteq T \subsetneq T' \\
         |S| \equiv |T| \mod 2}} 2^{|T'|-1-|T|}
    \right] = \\
    \sum_{\emptyset \neq T' \subseteq \{1,\dots,n\}}\delta_{T'}
    \left[
     \sum_{\substack{S \subseteq T \subseteq \{1,\dots,n\} \\
         |S| \not\equiv |T| \mod 2}}  2^{|T'|-1}
     - \sum_{\substack{S \subseteq T \subsetneq T' \\
         |S| \not\equiv |T| \mod 2}} 2^{|T'|-1-|T|}
    \right] +\delta_{T'}\frac{(-1)^{|T'|+|S|}-1}{2} \quad.
  \end{multline*}
  (Note that the last summand is equal to $- \delta_{T'}$
  for $|T'| \not\equiv |S| \mod 2$ and equal to $0$ otherwise.)

  Comparing coefficients for each $\delta_{T'}$ we obtain
  \begin{multline*}
    2^{n-|S|-1} \cdot 2^{|T'|-1}
    - \sum_{\substack{S \subseteq T \subsetneq T' \\
        |S| \equiv |T| \mod 2}} 2^{|T'|-1-|T|} = \\
    2^{n-|S|-1} \cdot 2^{|T'|-1}
     - \left( \sum_{\substack{S \subseteq T \subsetneq T' \\
         |S| \not\equiv |T| \mod 2}} 2^{|T'|-1-|T|} \right)
     +\frac{(-1)^{|T'|+|S|}-1}{2},
  \end{multline*}
  but this is true since
  \begin{equation*}
    (-1)^{|T'|-|S|} = (1-2)^{|T'|-|S|} = \sum_{\substack{S \subseteq T
\subseteq T' \\
        |T'| \equiv |T| \mod 2}} 2^{|T'|-|T|}
    -  \sum_{\substack{S \subseteq T \subseteq T' \\
         |T'| \not\equiv |T| \mod 2}} 2^{|T'|-|T|} \quad.
  \end{equation*}

  Suppose that $\bar{h}_{T}$ (for $\emptyset \neq T \subseteq
  \{1,\dots,n\}$) satisfy~\eqref{eq:existence} with equality, too.
  For $S \neq \emptyset$ inequality~\eqref{eq:existence} gives
  $$ \bar{h}_{S} = h_{S} + \bar{h}_{\{1,\dots,n\}} - h_{\{1,\dots,n\}}
  .$$
  But for $S=\emptyset$ inequality \eqref{eq:existence} yields
  $\bar{h}_{\{1,\dots,n\}} = h_{\{1,\dots,n\}}$ and therefore
  $\bar{h}_{T}=h_{T}$ for all $\emptyset \neq T \subseteq
  \{1,\dots,n\}$.

  This proves that~\eqref{eq:hT} is the only solution
  of~\eqref{eq:existence} that satisfies all inequalities with
  equality.

  Thus we find
  \begin{equation*}
      m \ge h_{\{1,\dots,n\}}
      \ge \sum_{\emptyset \neq T' \subseteq
        \{1,\dots,n\}}\delta_{T'}2^{|T'|-1}
      \ge \sum_{\emptyset \neq T' \subseteq
        \{1,\dots,n\}} 2^{|T'|-1}
      = \frac{1}{2}(3^{n}-1) \quad.
  \end{equation*}
\end{proof}

Next we prove a trade-off between the contrast of the
different images.

\begin{theorem}
  For $\emptyset \neq T \subseteq \{1,\dots,n\}$ let $\alpha_{T} =
 \frac{h_{T}-l_{T}}{m}$ be the contrast of the image $I_{T}$.
 The contrast levels of the images satisfy
 \begin{equation} \label{trade-off}
 \sum_{\emptyset \neq T \subseteq\{1,\dots,n\}} 2^{|T|-1}\alpha_{T}
 \le 1 \quad.
 \end{equation}

 Further let $\alpha'_{T}\ge 0$ (for $\emptyset \neq T \subseteq
 \{1,\dots,n\}$)
 satisfy \eqref{trade-off} .
 Then for every $\varepsilon > 0$ there exists a generalized visual
cryptography scheme
 with contrast levels
 $\alpha_{T}$ (for ${\emptyset \neq T \subseteq \{1,\dots,n\}}$) where
 $|\alpha_T - \alpha'_T| < \varepsilon$ for nonempty subsets $T$ of
 $\{1,\dots,n\}$.
\end{theorem}
\begin{proof}
  Let $\delta_{T} = h_{T}-l_{T}$. By~\eqref{eq:hT} we conclude
  $$m \ge h_{\{1,\dots,n\}} \ge \sum_{\emptyset \neq T \subseteq
    \{1,\dots,n\}} \delta_{T}2^{|T|-1}$$
  and therefore
  \begin{equation*}
      \sum_{\emptyset \neq T \subseteq\{1,\dots,n\}}
      2^{|T|-1}\alpha_{T} =  \frac{1}{m}\sum_{\emptyset \neq T
        \subseteq\{1,\dots,n\}} \delta_{T}2^{|T|-1} \\
      \le 1 \qquad.
  \end{equation*}

  Now assume \eqref{trade-off} holds for $\alpha'_{T}$.
  Then we choose $\delta_{T} \in \N$ and $M \in \N$ with
  $$ 0 \le \alpha'_{T} - \frac{\delta_{T}}{M} \le \epsilon \quad.$$
  By~\eqref{eq:hT} we know that there exists
  an extended visual cryptography scheme with
  contrast levels $h_{T}-l_{T}=\delta_{T}$ and minimal pixel
  expansion
  $$ m =\sum_{\emptyset \neq T \subseteq\{1,\dots,n\}}
  2^{|T|-1}\delta_{T}.$$

  Since $\frac{\delta_{T}}{M} \le \alpha'_{T}$ and $\alpha_{T'}$ satisfy
\eqref{trade-off}
  we find $m<M$.

  If we add useless subpixels (e.g. subpixels that are always black)
  to the extended visual cryptography scheme constructed above, we
  obtain a scheme with contrast
  $\alpha_{T}=\frac{\delta_{T}}{M}$. This proves the theorem.
\end{proof}


\section{Pixel Expansion and Contrast for the $\mathfrak{S}$-Extended
Scheme}
\label{sec:S}

In the previous section we proved that Construction \ref{droste} is
optimal if $\mathfrak{S}=\mathcal{P}(\{1,\dots,n\}) \backslash
\{\emptyset\}$.

If we set $\delta_{T}=1$ for $T \in \mathfrak{S}$ and
$\delta_{T}=0$ for $T \notin \mathfrak{S}$ in equation~\eqref{eq:hT}
we obtain the same contrast values as in Construction
\ref{droste}.
In this sense Construction \ref{droste}
can be viewed as a
$(\mathcal{P}(\{1,\dots,n\})\backslash\{\emptyset\})$-extended visual
cryptography scheme with degenerated contrast values.

\begin{theorem}\label{optimal}
  Let $\mathfrak{S} = \mathcal{P}(\{1,\dots,n\})\backslash
  \{\emptyset,\{1,\dots,n\}\}$ then Construction \ref{droste}
  is optimal if and only if $n$ is odd.
\end{theorem}
\begin{proof}
Let $\delta_T = 1$ for $\emptyset \neq T \subsetneq \{1 ,
\dots , n\}$ and $\delta_{\{1 , \dots , n\}} = 0$.
Then $h_T$ and $l_T = h_T - \delta_T$ as in (\ref{eq:hT}) satisfy
(\ref{eq:existence}) and these are the solutions given by
\ref{droste}.

Let us assume $n$ is even. We show that in this case
we can find a better construction in the sense that fewer subpixels
are required.
Setting $\overline{h}_T = h_T - 1$ and $\overline{l}_T = l_T - 1$
for $\emptyset \neq T \subseteq \{1,\dots,n\}$ and
$\overline{h}_{\emptyset}=\overline{l}_{\emptyset}=0$
we observe that inequality (\ref{eq:existence})
still holds for all $S \neq \emptyset$.
Thus the solution $x = (x_T)_{\emptyset \neq
T \subseteq \{1 , \dots ,n\}}$
of equation (\ref{matrix}) satisfies $x_T \geq 0$ for
$T \neq \{1 , \dots , n\}$.

The value of $x_{\{1 , \dots , n\}}$ will be non-negative unless
$r_{T}=\overline{h}_{T}$ for $|T|$ even and $r_{T}=\overline{l}_{T}$
for $|T|$ odd. In this special case equation~\eqref{eq:X-sol} gives the
solution
$x_{\{1,\dots,n\}}=-1$. To obtain a solution with positive
$x_{\{1,\dots,n\}}$ we adjust the value of $r_{\{1,\dots,n\}}$ from
$\overline{h}_{\{1,\dots,n\}}$ to $\overline{h}_{\{1,\dots,n\}}-1$.
(This is possible, since $\{1 , \dots , n\}
\not\in \mathfrak{S}$ which means that the number of black subpixels
in the stack of all transparencies does not matter.)
Now~\eqref{eq:X-sol} reveals the solution
\begin{align*}
  x_{\{1,\dots,n\}} &=
    \sum_{T \subseteq \{1,\dots,n\}} (-1)^{|T|+1} r_{T}  \\
    &=\sum_{\substack{T \subseteq \{1,\dots,n\} \\ |T| \text{ odd}}}
    \overline{l}_{T} -
    \left(\sum_{\substack{T \subseteq \{1,\dots,n\} \\ |T| \text{ even}}}
    \overline{h}_{T}\right) + 1 = -1+1=0. \\
    \intertext{For $S \neq \{1,\dots,n\}$ and $|S|$ even we obtain}
  x_{S} &= \sum_{\{1,\dots,n\}\backslash S \subseteq T \subseteq
    \{1,\dots,n\}} (-1)^{|T|+|S|+n+1} r_{T} \\
  &=\sum_{\substack{\{1,\dots,n\}\backslash S \subseteq T \subseteq
      \{1,\dots,n\} \\ |T| \text{ odd}}}
  \overline{l}_{T} -
  \left(\sum_{\substack{\{1,\dots,n\}\backslash S \subseteq T \subseteq
      \{1,\dots,n\} \\ |T| \text{ even}}}
  \overline{h}_{T}\right) + 1 = 0+1=1.
  \intertext{For $|S|$ odd we obtain}
  x_{S} &= \sum_{\{1,\dots,n\}\backslash S \subseteq T \subseteq
    \{1,\dots,n\}} (-1)^{|T|+|S|+n+1} r_{T} \\
  &=\left(\sum_{\substack{\{1,\dots,n\}\backslash S \subseteq T \subseteq
      \{1,\dots,n\} \\ |T| \text{ even}}}
  \overline{h}_{T}\right) - 1 -
  \sum_{\substack{\{1,\dots,n\}\backslash S \subseteq T \subseteq
      \{1,\dots,n\} \\ |T| \text{ odd}}}
  \overline{l}_{T}  \\
  &> \sum_{\substack{\{1,\dots,n\}\backslash S \subseteq T \subseteq
      \{1,\dots,n\} \\ |T| \text{ even}}}
  \overline{l}_{T} -
  \left(\sum_{\substack{\{1,\dots,n\}\backslash S \subseteq T \subseteq
      \{1,\dots,n\} \\ |T| \text{ odd}}}
  \overline{h}_{T}\right) -1 = -1.\\
\end{align*}
Thus all possible values of $r$ lead to non-negative solutions for
$x$, hence an $\mathfrak{S}$-extended visual cryptography
scheme with $m = \overline{h}_{\{1,\dots,n\}}<h_{\{1,\dots,n\}}$
exists, i.e. the solution given by \ref{droste} is not optimal.

Now we assume $n$ is odd. Let $\overline{h}_{T}$ and
$\overline{l}_{T}$ ($\emptyset \neq T \subsetneq \{1 , \dots ,n\})$)
be a solution of \eqref{eq:existence} different from
$h_T$ and $l_T$. Let $\overline{h}_{\{1 , \dots ,n\}}$ be the maximal
number of black subpixels occurring when all the transparencies are
stacked together.

Since $\overline{h}_{T}$ and
$\overline{l}_{T}$ solve inequality (\ref{eq:existence}), we can apply
the same arguments leading from (\ref{eq:existence}) to \eqref{eq:hT}
and find
\begin{equation}\label{eq:h-diff}
\overline{h}_{\{1 , \dots ,n\}} - \overline{h}_{T} \geq
h_{\{1 , \dots ,n\}} - h_T
\end{equation}
for all $\emptyset \neq T \subseteq \{1 , \dots ,n\}$. But for $S =
\emptyset$ in (\ref{eq:existence}), we obtain
$$\sum_{\substack{T \subseteq \{1,\dots,n\} \\ |T| \text{ even}}}
    \overline{h}_{T} \leq \sum_{\substack{T \subseteq \{1,\dots,n\} \\
        |T| \text{ odd}}} \overline{l}_{T}$$
and hence, together with \eqref{eq:h-diff},
$$\overline{l}_{\{1 , \dots ,n\}} \geq l_{\{1 , \dots ,n\}} =
h_{\{1 , \dots ,n\}}.$$
But by definition $\overline{h}_{\{1 , \dots ,n\}} \geq
\overline{l}_{\{1 , \dots ,n\}}$ and therefore this solution needs at
least $h_{\{1 , \dots ,n\}}$ subpixels, i.e. the solution given by
\ref{droste} is optimal.
\end{proof}

 We notice that in the first part of the proof of \ref{optimal}
 we did not need the assumption $\delta_T = 1$ for $\emptyset \neq T
 \subsetneq \{1, \dots , n\}$. In fact, it is sufficient to assume
 $h_T \neq 0$ where $h_T$ is defined by (\ref{eq:hT}). A short
 calculation proves that this is the case if $\mathfrak{S}
 \not\subseteq \mathcal{P}(S)$
 for a proper subset $S$ of $\{1,\dots,n\}$.
 Thus the first part of \ref{optimal} proves:

\begin{corollary}\label{optimal2}
 For $\mathfrak{S}
 \subseteq \mathcal{P}(\{1,\dots,n\}) \backslash \{\emptyset,
 \{1,\dots,n\} \}$, $\mathfrak{S} \not\subseteq S$ for any proper subset
 $S$ of $\{1,\dots,n\}$ and $n$ even, Construction \ref{droste} is not
 optimal.
\end{corollary}

Even more general we can prove:

\begin{theorem} \label{optimal2.5}
  Let $\mathfrak{S} \subseteq \mathcal{P}(\{1,\dots,n\}) \backslash
  \{\emptyset\}$. Let us assume that there exists a nonempty subset
  $T \in \mathcal{P}(\{1,\dots,n\})
  \backslash \mathfrak{S}$ with $|T|$ even and that $\mathfrak{S} \cap
  \mathcal{P}(T) \not \subseteq \mathcal{P}(T')$ for each proper subset
  $T'$ of $T$. Then Construction \ref{droste} is not optimal.
\end{theorem}

\begin{proof}
We apply \ref{optimal2} to $\mathfrak{S} \cap \mathcal{P}(T)$. For a
proper subset $T'$ of $T$ we define $\overline{h}_{T'}$ and
$\overline{l}_{T'}$ as in the proof of \ref{optimal}. Formally we
define $\overline{h}_T = h_T - 2$ and  $\overline{l}_T = h_T -1$ where
$h_T$ is defined by \eqref{eq:hT}. Let $\delta_T = -1$ and for $T \neq
S \subseteq \{1,\dots,n\}$ we define $\delta_S = 0$ for $S \not\in
\mathfrak{S}$ and $\delta_S = 1$ otherwise.

Corresponding to \eqref{eq:hT} we define
\begin{equation*}
\hat{h}_S =
\begin{cases}
\overline{l}_T + \sum\limits_{\substack{\emptyset \neq S' \subsetneq
    \{1,\dots,n\} \\ S' \not\subseteq T}} \delta_{S'}2^{|S'|-1} -
    \sum\limits_{S \subsetneq S' \subseteq \{1,\dots,n\}}
    \delta_{S'}2^{|S'|-1-|S|} &
    \text{for $S \not\subseteq T$} \\
\\
\overline{h}_S + \sum\limits_{\substack{\emptyset \neq S' \subsetneq
    \{1,\dots,n\} \\ S' \not\subseteq T}} \delta_{S'}2^{|S'|-1} -
    \sum\limits_{\substack{S \subsetneq S' \subseteq \{1,\dots,n\}\\S'
    \not\subseteq T}} \delta_{S'}2^{|S'|-1-|S|} &
    \text{for $S \subseteq T$} \\
\end{cases}
\end{equation*}
for $\emptyset \neq S \subseteq \{1, \dots , n\}$ and $\hat{l}_S =
\hat{h}_S - \delta_S$.

An easy but tedious
calculation corresponding to the one that was used in the proof of
Theorem \ref{drostegut} shows that $\hat{h}_S, \hat{l}_S$ satisfy
\eqref{eq:existence}. Hence an $\mathfrak{S}$-extended visual
cryptography system with contrast levels $\hat{h}_S, \hat{l}_S$
exists. This scheme has a smaller pixel-expansion than the scheme
given by Construction \ref{droste}.
\end{proof}

\section{Conclusions and further remarks}

Equation~\eqref{eq:X-sol} gives us a simple method to construct an
$\mathfrak{S}$-extended visual cryptography scheme with given
contrast values $l_{T}$ and $h_{T}$. Furthermore, for fixed $n$ and
$\mathfrak{S}$, equation~\eqref{eq:X-sol} leads to a linear programming
problem which describes all possible $\mathfrak{S}$-extended visual
cryptography schemes. For small values of $n$ this problem can easily
be solved.

In this article we have given a full solution for the special cases
$\mathfrak{S}=\mathcal{P}(\{1,\dots,n\}) \backslash \{\emptyset\}$
and $\mathfrak{S}=\mathcal{P}(\{1,\dots,n\}) \backslash \{\emptyset,
\{1,\dots,n\}\}$. We close by presenting the following open problems:

\begin{enumerate}
\item We conjecture:
 \begin{quote}
  Let $\mathfrak{S} \subseteq \mathcal{P}(\{1,\dots,n\})\backslash
  \{\emptyset\}$.
  Then
  Construction \ref{droste} is optimal for an $\mathfrak{S}$-extended
  visual cryptography scheme if and only if for all $\emptyset \neq T
  \not\in \mathfrak{S}$ we have either $|T|$ odd or $\mathfrak{S} \cup
  \mathcal{P}(T) \subseteq \mathcal{P}(T')$ for some proper subset $T'$
  of $T$.
 \end{quote}
\item An even harder problem is a full characterization of
  $\mathfrak{S}$-extended visual cryptography schemes with minimal
  pixel expansion for arbitrary
  subsets $\mathfrak{S}$ of $\mathcal{P}(\{1,\dots,n\})$, i.e. to find
  a formula for the minimal pixel expansion depending on
  $\mathfrak{S}$.
\end{enumerate}

\scriptsize

Andreas Klein\\
Universit\"at Kassel\\
Fachbereich 17 (Mathematik und Informatik)\\
D-34109 Kassel\\
{\tt klein@mathematik.uni-kassel.de}

Markus Wessler\\
Universit\"at Kassel\\
Fachbereich 17 (Mathematik und Informatik)\\
D-34109 Kassel\\
{\tt wessler@mathematik.uni-kassel.de}

\end{document}